\documentclass[11pt]{amsart}
\usepackage[utf8]{inputenc}
\usepackage[usenames,dvipsnames]{xcolor}
\usepackage[draft]{pgf}
\usepackage{listings}
\usepackage{float}
\usepackage{caption}

\usepackage{pgfplots}
\pgfplotsset{compat=1.15}
\usepackage{mathrsfs}
\usetikzlibrary{arrows}

\usepackage[margin=1.2in, centering]{geometry}
\usepackage{amsmath,amssymb,latexsym,mathtools}
\usepackage{dsfont}
\usepackage{nccmath}

\usepackage{amsthm}
\usepackage{amssymb}
\newtheorem{thm}{Theorem}[section]

\numberwithin{equation}{section}

\theoremstyle{definition}

\newtheorem*{acknowledgement}{Acknowledgements}

\usepackage{graphicx}
\DeclareGraphicsExtensions{.pdf,.png,.jpg}
\usepackage{subfig}

\counterwithin*{equation}{section}

\usepackage{hyperref}
\hypersetup{linktocpage}

\begin{document}

\definecolor{qqttcc}{rgb}{0.,0.2,0.8}
	\title{Near optimal thresholds for existence of dilated configurations in $\mathbb{F}_q^d$}
	\author{Pablo Bhowmik and Firdavs Rakhmonov}
	
	\date{\today}

	
        \address{Department of Mathematics, University of Rochester, Rochester, NY, USA}
	\email{pbhowmik@ur.rochester.edu}

        \address{Department of Mathematics, University of Rochester, Rochester, NY, USA}
	\email{frakhmon@ur.rochester.edu}
 
	\keywords{}
	\maketitle

\begin{abstract} 
    Let $E\subset\mathbb{F}_q^d$ and $\lVert \cdot \rVert:\mathbb{F}_q^d\to \mathbb{F}_q$ defined as $\lVert \alpha\rVert\coloneqq \alpha_1^2+\dots+\alpha_d^2$ if $\alpha=(\alpha_1,\dots,\alpha_d)\in \mathbb{F}_q^d$, where $\mathbb{F}_q^d$ is the $d$-dimensional vector space over the finite field $\mathbb{F}_q$ with $q$ elements. Let $k\geq 1$ and $A$ is a nonempty subset of $\{(i,j):1\leq i<j\leq k+1\}$. In this paper, we prove that for any nonzero square element $r\in \mathbb{F}_q$ and $\lvert E\rvert\gg_kq^{d/2}$ one can find two $(k+1)$-tuples in $E$ such that one of them is dilated by $r$ with respect to the other only on $|A|$ edges. More precisely, there exist $(x_1,\dots,x_{k+1})\in E^{k+1}$ and $(y_1,\dots,y_{k+1})\in E^{k+1}$ such that $\lVert y_i-y_j\rVert=r\lVert x_i-x_j\rVert$ if $(i,j)\in A$ and $x_i\neq x_j, y_i\neq y_j$ if $1\leq i<j\leq k+1$. In particular, we present two proofs first of which utilizes the machinery of group actions which is typically used to handle problems of this nature, whereas the second one is more combinatorial in nature and is based on the averaging argument. Moreover, we show that in two dimensions the threshold $d/2$ is sharp when $q\equiv 3 \pmod 4$.  
    
    \smallskip

    As a corollary of this result, varying the underlying set $A$ we obtain thresholds for existence of dilated $k$-paths, $k$-cycles, and $k$-stars $(k\geq 3)$. These results partly generalize some results of the second author. 
\end{abstract} 

\bigskip

\section{Introduction}

Let $\mathbb F_q^d, d\ge 2,$ be the $d$-dimensional vector space over the finite field $\mathbb F_q$ with $q$ elements. Throughout this paper, we assume that $q$ is a power of odd prime $p.$ 
Given a set $E$ in $\mathbb F_q^d$,  the distance set $\Delta(E)$ is defined by
$$ \Delta(E):=\{\lVert x-y \rVert \in \mathbb F_q: x,y \in E\},$$
where $||\alpha||=\sum_{i=1}^d \alpha_i^2$ for $\alpha=(\alpha_1, \ldots, \alpha_d)\in \mathbb F_q^d.$ 

\smallskip

In the finite field setting, the Falconer distance problem asks for the smallest exponent $\alpha>0$ such that  for any $E\subset\mathbb F_q^d$ with $|E|\ge C q^\alpha,$  we have $ |\Delta(E)|\ge c q$, where $C>1$ denotes a sufficiently large constant, and  $0<c\leq 1$ denotes some constant independent of $q$ and $|E|,$ where $|E|$ denotes the cardinality of $E.$  This problem was proposed by Iosevich and Rudnev \cite{MR2336319} as a finite field analogue of the Falconer distance problem in the Euclidean space. We notice that the formulation of the finite field Falconer problem was also motivated on the Erd\H{o}s distinct distances problem over finite fields due to  Bourgain, Katz and Tao \cite{MR2053599}, and that is why the problem is also called the Erd\H{o}s-Falconer distance problem.  
We refer readers to \cite{MR4201782, MR4297185, MR3961084, MR15796, MR834490, MR4055179, MR2721878, MR3272924, MR2399013} for precise definitions, background knowledge, and recent progress on the Erd\H{o}s distinct distances problem and  the Falconer distance problem in the Euclidean setting.  

\smallskip

One can consider a strong version of the Erd\H{o}s-Falconer distance problem which is called the Mattila-Sj\"olin distance problem over finite fields. It asks for the smallest threshold $\beta >0$ such that for any $E\subset\mathbb F_q^d$ with $|E|\ge Cq^\beta$ we have $ \Delta(E)=\mathbb F_q.$ 
Using Fourier analytic technique and the Kloosterman sum estimates, Iosevich and Rudnev  obtained the threshold $(d+1)/2$ for all dimensions $d\ge 2.$

\begin{thm} [Iosevich and Rudnev, \cite{MR2336319}]
\label{IRthm} 
If $E\subset \mathbb F_q^d \ (d\ge 2)$ 
and  $|E|>2 q^{(d+1)/2}$, then $\Delta(E)=\mathbb F_q.$
\end{thm}

The threshold $(d+1)/2$ in Theorem \ref{IRthm} is the best currently known result on the Mattila-Sj\"olin distance problem over finite fields for all dimensions $d\ge 2.$ It is considered as a challenging problem to improve the $(d+1)/2$ threshold. We notice that in odd dimensions, it gives the optimal threshold, which was proven in \cite{MR2775806}.

\smallskip
 
However, in even dimensions,  it has been believed that the exponent $(d+1)/2$ can be improved but reasonable evidence or conjecture has not been stated in literature. In two dimensions, Murphy and Petridis \cite{MR3955393} showed that the threshold  cannot be lower than $4/3$ for the Mattila-Sj\"olin distance problem over finite fields.

\smallskip

On the other hand, Iosevich and Rudnev \cite{MR2336319} conjectured that the threshold $(d+1)/2$ can be lower to $d/2$ for the Erd\H{o}s-Falconer distance problem in even dimensions, and  in two dimensions the threshold $4/3$ was proven by the authors in \cite{MR2917133} (see also \cite{MR3592595}). Moreover, if $q$ is a prime number, then the exponent $4/3$ was improved to $5/4$ by Murphy,  Petridis,  Pham,  Rudnev, and  Stevenson \cite{MR4411329}. We also notice that the threshold $(d+1)/2$ cannot be improved for the Erd\H{o}s-Falconer distance problem in odd dimensions (see also \cite{MR2775806}). 

\smallskip

The distance problems over finite fields have been extended  in various directions. Although lots of variants of the distance problems were extensively studied, the threshold $d/2$ for the set $E$ in $\mathbb F_q^d$ had not been addressed for any distance type problems until Iosevich, Koh and Parshall \cite{MR3959878} studied the Mattila-Sj\H{o}lin problem for the quotient set of the distance set over finite fields. 

\smallskip

\textbf{The Mattila-Sj\H{o}lin problem for the quotient set of the distance set.}\label{prIKP} If $E$ is a subset of $\mathbb{F}_q^d \ (d\geq 2)$, then the quotient set of the distance set, denoted by $\frac{\Delta(E)}{\Delta(E)},$ is defined as follows:
$$\frac{\Delta(E)}{\Delta(E)}:=\left\{ \frac{\lVert x-y \rVert}{\lVert z-w \rVert}: x,y, z, w\in E, ~ \lVert z-w \rVert\ne 0\right\}.$$

So one can ask the following question: determine the smallest exponent $\gamma>0$ such that for any set $E \subset \mathbb F_q^d$ with $|E|\ge C q^\gamma$, we have
\begin{equation}\label{Qeq} \frac{\Delta(E)}{\Delta(E)}=\mathbb F_q.\end{equation}

The aforementioned authors obtained the threshold $d/2$ in even dimensions on the Mattila-Sj\H{o}lin problem for the quotient set of the distance set over finite fields. More precisely, they proved the following result.

\begin{thm} [Theorems 1.1, 1.2, \cite{MR3959878}]
\label{IKPmthm} 
If $E \subset \mathbb F_q^d \ (d\ge 2)$, then the following statements hold:
\begin{enumerate} \item  [(i)]If $d\ge 2$ is even and $ |E|\ge 9 q^{d/2},$ then  $\dfrac{\Delta(E)}{\Delta(E)}=\mathbb F_q.$

\item [(ii)] If $d\ge 3$ is odd and $|E|\ge 6 q^{d/2},$ then  $\dfrac{\Delta(E)}{\Delta(E)}\supseteq(\mathbb F_q)^2,$
where $(\mathbb F_q)^2 := \{a^2: a\in \mathbb F_q\}.$  
\end{enumerate}
\end{thm}

We notice that recently Theorem \ref{IKPmthm} has been extended and generalized with improved constants to the general non-degenerate quadratic distances by Iosevich, Koh and the second author (see \cite{IKR2023}).

\smallskip

We shall write $\mathbb F_q^*$ for the set of all non-zero elements in $\mathbb F_q.$ 

\smallskip

Equality ${\Delta(E)}/{\Delta(E)}=\mathbb{F}_q$ implies that for any $r\in \mathbb{F}_q^{*}$ one can find $(x,y)\in E^2$ and $(x',y')\in E^2$ such that $\lVert x-y\rVert\neq 0$ and $\lVert x'-y' \rVert=r \lVert x-y \rVert$. In other words, if $r\in \mathbb{F}_q^*$ and $E\subset \mathbb{F}_q^d$ with $|E|\geq 9q^{d/2}$, then one can find a pair of edges in the complete graph $K_{|E|}$ with vertex set $E$ such that one of them is dilated by $r$ with respect to the other.

\smallskip

A natural question arises whether it is possible to generalize this result to other subgraphs of the complete graph $K_{|E|}$ with vertex set $E$. The case of 2-paths, 4-cycles and simplexes have been studied by the second author (see \cite{Rakh2022}).

\smallskip

Let $r\in \mathbb{F}_q^{*}$, $E\subset \mathbb{F}_q^d$ and assume that $A$ is a nonempty subset of $\{(i,j):1\leq i<j\leq k+1\}$, where $k\geq 1$. One can ask how large $E$ needs to be to guarantee the existence of two point-configurations $(x_1,\dots, x_{k+1})\in E^{k+1}$ and $(y_1,\dots,y_{k+1})\in E^{k+1}$ such that $x_i\neq x_j,\ y_i\neq y_j$ if $1\leq i< j\leq k+1$ and $\lVert y_i-y_j\rVert=r\lVert x_i-x_j\rVert$ if $(i,j)\in A$.

\smallskip

\begin{figure}[htbp]
\centering
\usetikzlibrary{arrows}
\definecolor{qqqqff}{rgb}{0.,0.,1.}
\definecolor{ududff}{rgb}{0.30196078431372547,0.30196078431372547,1.}
\begin{tikzpicture}[scale=0.75][line cap=round,line join=round,>=triangle 45,x=1.0cm,y=1.0cm]
\clip(-0.5,-3.4) rectangle (17.,-0.5);
\draw [line width=1.pt,color=qqqqff,fill=qqqqff,fill opacity=0.0] (1.84,-2.96)-- (2.96,-1.18);
\draw [line width=1.pt,color=qqqqff,fill=qqqqff,fill opacity=0.0] (2.96,-1.18)-- (5.7,-1.);
\draw [line width=1.pt,color=qqqqff,fill=qqqqff,fill opacity=0.0] (5.7,-1.)-- (5.78,-3.18);
\draw [line width=1.pt,color=qqqqff,fill=qqqqff,fill opacity=0.0] (5.78,-3.18)-- (1.84,-2.96);
\draw [line width=1.pt,color=qqqqff,fill=qqqqff,fill opacity=0.0] (1.84,-2.96)-- (5.7,-1.);
\draw [line width=1.pt,color=qqqqff,fill=qqqqff,fill opacity=0.0] (8.68,-3.08)-- (14.56,-3.18);
\draw [line width=1.pt,color=qqqqff,fill=qqqqff,fill opacity=0.0] (14.56,-3.18)-- (13.22,-0.76);
\draw [line width=1.pt,color=qqqqff,fill=qqqqff,fill opacity=0.0] (13.22,-0.76)-- (10.,-0.78);
\draw [line width=1.pt,color=qqqqff,fill=qqqqff,fill opacity=0.0] (10.,-0.78)-- (8.68,-3.08);
\draw [line width=1.pt,color=qqqqff,fill=qqqqff,fill opacity=0.0] (8.68,-3.08)-- (13.22,-0.76);
\draw [fill=ududff] (1.84,-2.96) circle (2.5pt);
\draw[color=black] (1.45,-3.08) node {\scalebox{1}{$x_1$}};
\draw [fill=ududff] (2.96,-1.18) circle (2.5pt);
\draw[color=black] (2.55,-1) node {\scalebox{1}{$x_2$}};
\draw [fill=ududff] (5.7,-1.) circle (2.5pt);
\draw[color=black] (6.2,-0.92) node {\scalebox{1}{$x_3$}};
\draw [fill=ududff] (5.78,-3.18) circle (2.5pt);
\draw[color=black] (6.2,-3.2) node {\scalebox{1}{$x_4$}};
\draw [fill=ududff] (8.68,-3.08) circle (2.5pt);
\draw[color=black] (8.3,-3.2) node {\scalebox{1}{$y_1$}};
\draw [fill=ududff] (14.56,-3.18) circle (2.5pt);
\draw[color=black] (15,-3.2) node {\scalebox{1}{$y_4$}};
\draw [fill=ududff] (13.22,-0.76) circle (2.5pt);
\draw[color=black] (13.7,-0.7) node {\scalebox{1}{$y_3$}};
\draw [fill=ududff] (10.,-0.78) circle (2.5pt);
\draw[color=black] (9.6,-0.7) node {\scalebox{1}{$y_2$}};
\end{tikzpicture}
\label{2-chains}
\caption{Pair of quadrilaterals with dilation ratio $r\in \mathbb{F}_p^{*}$. }
\label{2chainsfigure}
\end{figure}
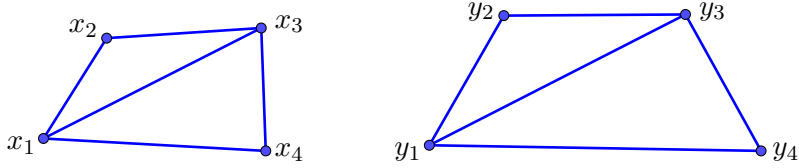

For instance, let $r\in \mathbb{F}_p^{*}$ and $E\subset \mathbb{F}_p^2$. We would like to determine the smallest exponent $\gamma>0$ such that for any $E\subset \mathbb{F}_p^2$ with $|E|\geq Cp^{\gamma}$ there are 4-tuples $(x_1,x_2,x_3,x_4)\in E^4$ and $(y_1,y_2,y_3,y_4)\in E^4$ such that $x_i\neq x_j,\ y_i\neq y_j$ if $1\leq i< j\leq 4$ and $\lVert y_i-y_j\rVert=r\lVert x_i-x_j \rVert$ if $(i,j)\in A$, where $A=\{(1,2),(2,3),(3,4),(1,4),(1,3)\}$ (see Figure \ref{2chainsfigure}).

One of the most remarkable results in this topic is due to Bourgain (see \cite{Bourgain86}). It states that if $A$ is a set of positive upper Lebesgue density and $u^1, \ldots, u^k$ are vectors in $\mathbb{R}^d \ (k\le d)$, then there exists $r_0>0$ such that for all $r>r_0$ and any $\delta>0$, there exist $x^1, \ldots, x^{k+1}$ in the $\delta$-neighborhood of $A$ forming a simplex similar to $\{0, u^1,\ldots, u^k\}$ with the scaling coefficient  $r$. 

\smallskip

In this paper, we study the discrete setting of this type problem, namely, in the vector space over finite fields. 

\smallskip

We solve this problem for $(k+1)$-point configurations $(k\geq 1)$ in $\mathbb{F}_q^d$ with fixed nonempty subset $A$ of $\{(i,j): 1\leq i<j\leq k+1\}$. 
We provide two proofs and the first one utilizes the machinery of group actions and boils down to the investigation of $L^{k}$ and $L^{k+1}$-norms of the following counting function: 
\begin{equation*}
    \lambda_{r,\theta}(z)\coloneqq \big\{(u,v)\in E^2: u-\sqrt{r}\theta v=z\big\},    
\end{equation*}
 where $r\in (\mathbb{F}_q)^2\setminus \{0\}$, $\theta\in \mathrm{O}_d(\mathbb{F}_q)$ and $z\in \mathbb{F}_q^d$. Meanwhile, the second proof is more combinatorial in nature and is based on the averaging argument.

\smallskip

Let us formulate the main result of this paper.

\begin{thm}
\label{mainthm}
    Suppose $r\in (\mathbb{F}_q)^2\setminus \{0\}$ and $\varnothing \neq A\subset \big\{(i,j):1\leq i<j\leq k+1\big\}$, where $k\geq 1$. If $E\subset \mathbb{F}_q^d$ with $\lvert E\rvert\geq 2kq^{d/2}$, then there exist $(x_1,\dots,x_{k+1})\in E^{k+1}$ and $(y_1,\dots,y_{k+1})\in E^{k+1}$ such that $\lVert y_i-y_j \rVert=r\lVert x_i-x_j \rVert$ if $(i,j)\in A$ and $x_i\neq x_j, y_i\neq y_j$ if $1\leq i<j\leq k+1$.
\end{thm}

\smallskip

As a simple corollary, varying the underlying set $A$ we obtain thresholds for existence of dilated $k$-paths, $k$-cycles and $k$-stars $(k\geq 3)$. For example, the result about $k$-paths $(k\geq 3)$ generalizes Theorem 1.5 in \cite{Rakh2022}.

\smallskip

The rest of this paper will be organized to prove Theorem \ref{mainthm}.

\begin{acknowledgement}
    Authors would like to thank Kaave Hosseini and Alexander Iosevich for suggesting the interesting problem and for helpful discussions.
\end{acknowledgement}  

\section{Notations}
\label{sec:2}

We recall some basic notations which we will use throughout the paper.

\medskip

Let $\mathrm{O}_d(\mathbb{F}_q)$ denote the group of orthogonal $d\times d$ matrices with entries in $\mathbb{F}_q$.

\smallskip

Let $\mathbb{F}_q^*$ denote the set of nonzero elements in $\mathbb{F}_q$, i.e., $\mathbb{F}^*_q=\big\{a\in \mathbb{F}_q: a\neq 0 \big\}$.

\smallskip

Let $(\mathbb{F}_q)^2$ denote the set of quadratic residues in $\mathbb{F}_q$, i.e., $(\mathbb{F}_q)^2=\big\{a^2: a\in \mathbb{F}_q \big\}$.

\smallskip 
 
Consider the map $\lVert \cdot\rVert:\mathbb{F}_q^d\to \mathbb{F}_q$ which is defined by $\lVert \alpha \rVert\coloneqq \alpha_1^2+\dots+\alpha_d^2$, where $\alpha=(\alpha_1,\dots,\alpha_d)\in \mathbb{F}_q^d$. This mapping is not a norm since we do not impose any metric structure on $\mathbb{F}_q^d$ but it shares the important feature of the Euclidean norm: it is invariant under orthogonal matrices.

\smallskip

If $X$ is a finite set, then let $|X|$ denote the size (the cardinality) of $X$.

\smallskip

If $n\in \mathbb{N}$, then $[n]\coloneqq \{1,\dots,n\}$.

\smallskip

Let $\lfloor x\rfloor$ denote the greatest integer less than or equal to $x$.

\smallskip

We write $X\gg_{k} Y$ as shorthand for the inequality $X\geq C_{k}Y$ for some constant $C_{k}>0$ depending only on $k$.

\bigskip

\section{First Proof of Theorem \ref{mainthm}}
\label{sec:3}

Let $k\geq 1$ and $\varnothing \neq A\subset \big\{(i,j):1\leq i<j\leq k+1\big\}$. For $r\in (\mathbb{F}_q)^2\setminus \{0\}$, $z\in \mathbb{F}_q^d$ and $\theta \in \mathrm{O}_d(\mathbb{F}_q)$, we consider the following counting function:
\begin{equation}
\label{definitionofcountingfunction}
    \lambda_{r,\theta}(z)\coloneqq \left\lvert\big\{(u,v)\in E^2: u-\sqrt{r}\theta v=z\big\}\right\rvert.
\end{equation}

For $p\geq 1$, we define its $L^{p}$-norm as follows:
\begin{equation*}
    \lVert \lambda_{r,\theta}(z)\rVert_{p}^{p}\coloneqq \sum_{\theta\in \mathrm{O}_d(\mathbb{F}_q)}\sum_{z\in \mathbb{F}_q^d}\lambda_{r,\theta}(z)^{p}.
\end{equation*}

From \eqref{definitionofcountingfunction}, we see that
\begin{equation*}
    \lambda_{r,\theta}(z)^{k+1}=\big\lvert\big\{(u_1,\dots,u_{k+1},v_1,\dots,v_{k+1})\in E^{2k+2}: u_i-\sqrt{r}\theta v_i=z,\ i\in [k+1]\big\}\big\rvert.
\end{equation*}

\smallskip

Then, we obtain

\begin{equation*}
\begin{split}
   \lVert \lambda_{r,\theta}(z)\rVert_{k+1}^{k+1}
   =\sum_{\theta,z}\big\lvert\big\{(u_1,\dots,u_{k+1},v_1,\dots,v_{k+1})\in E^{2k+2}: u_i-\sqrt{r}\theta v_i=z,\ i\in [k+1]\big\}\big\rvert\\
   =\sum_{\theta}\big\lvert\big\{(u_1,\dots,u_{k+1},v_1,\dots,v_{k+1})\in E^{2k+2}: u_i-u_j=\sqrt{r}\theta (v_i-v_j),\ 1\leq i< j\leq k+1 \big\}\big\rvert.
\end{split}
\end{equation*}

Let us introduce the following notation:

\begin{equation*}
    \Lambda_{\theta}(r)\coloneqq \Big\{(u_1,\dots,u_{k+1},v_1,\dots,v_{k+1})\in E^{2k+2}: u_i-u_j=\sqrt{r}\theta (v_i-v_j),\ 1\leq i< j\leq k+1\Big\}.
\end{equation*}

\smallskip

Using the introduced notation one can write the $L^{k+1}$-norm of the function $\lambda_{r,\theta}(z)$ in terms of $|\Lambda_{\theta}(r)|$, i.e.,

\begin{equation*}
    \lVert \lambda_{r,\theta}(z)\rVert_{k+1}^{k+1}=\sum_{\theta\in \mathrm{O}_d(\mathbb{F}_q)}|\Lambda_{\theta}(r)|.
\end{equation*}

Let us consider the following auxiliary sets:

\begin{equation*}
     \mathcal{N}_{A,\theta}(r):=\left\{(u_1,\dots,u_{k+1},v_1,\dots,v_{k+1})\in E^{2k+2}: \begin{array}{l} u_i-u_{j}=\sqrt{r}\theta(v_i-v_j),\ (i,j)\in A,\\ v_i\neq v_j,\ 1\leq i< j\leq k+1\end{array}\right\},    
\end{equation*}

\smallskip

\begin{equation*}
     \mathcal{N}_{\theta}(r):=\left\{(u_1,\dots,u_{k+1},v_1,\dots,v_{k+1})\in E^{2k+2}: \begin{array}{l} u_i-u_{j}=\sqrt{r}\theta(v_i-v_j),\\ v_i\neq v_j,\ 1\leq i<j\leq k+1\end{array}\right\}.    
\end{equation*}

\smallskip

In order to prove Theorem \ref{mainthm} we need to show that $\lvert\mathcal{S}_A\rvert>0$, where

\begin{equation*}
        \mathcal{S}_A\coloneqq\left\{(x_1,\dots,x_{k+1},y_1,\dots,y_{k+1})\in E^{2k+2}: \begin{array}{l}\lVert y_i-y_{j}\rVert=r\lVert x_i-x_{j}\rVert,\ (i,j)\in A,\\
        x_i\neq x_j,\ y_i\neq y_j,\ 1\leq i<j\leq k+1 \end{array}\right\}.\\
\end{equation*}

\smallskip

It is easy to verify that for each $\theta\in \mathrm{O}_d(\mathbb{F}_q)$, we have $\mathcal{N}_{A,\theta}(r)\subset \mathcal{S}_A$ and hence 
\begin{equation}
\label{firstmaininequality}
    \lvert \mathcal{S}_A\rvert\geq \frac{1}{\lvert \mathrm{O}_d(\mathbb{F}_q)\rvert}\sum_{\theta\in \mathrm{O}_d(\mathbb{F}_q)}|\mathcal{N}_{A,\theta}(r)|.
\end{equation}

We notice that for each $\theta\in \mathrm{O}_d(\mathbb{F}_q)$, we have $\mathcal{N}_{\theta}(r)\subset \mathcal{N}_{A,\theta}(r)$. Comparing this observation with \eqref{firstmaininequality}, we obtain 
\begin{equation}
\label{secondmaininequality}
    \lvert \mathcal{S}_A\rvert\geq \frac{1}{\lvert \mathrm{O}_d(\mathbb{F}_q)\rvert}\sum_{\theta\in \mathrm{O}_d(\mathbb{F}_q)}|\mathcal{N}_{\theta}(r)|.
\end{equation}

For each pair $(\alpha,\beta)$ such that $1\leq \alpha<\beta\leq k+1$, we define the following set: 
\smallskip
\begin{equation*}
     A_{\alpha,\beta}:=\left\{(u_1,\dots,u_{k+1},v_1,\dots,v_{k+1})\in E^{2k+2}: \begin{array}{l} u_i-u_{j}=\sqrt{r}\theta(v_i-v_j),\\ 
1\leq i<j\leq k+1,\ v_{\alpha}=v_{\beta}\end{array}\right\}.
\end{equation*}

\smallskip

One can see that the following set equality holds:
\begin{equation*}
    \Lambda_{\theta}(r)\setminus \bigcup_{1\leq \alpha<\beta\leq k+1}A_{\alpha,\beta}=\mathcal{N}_{\theta}(r).
\end{equation*}

Hence, we obtain the following inequality:
\smallskip
\begin{equation}
\label{eq6.4}
    |\mathcal{N}_{\theta}(r)|\geq |\Lambda_{\theta}(r)|-\sum \limits_{1\leq \alpha<\beta\leq k+1}|A_{\alpha,\beta}|.
\end{equation}

One can verify that if $(\alpha,\beta)$ with $1\leq \alpha<\beta\leq k+1$, then we have 
\begin{equation}
\label{eq6.5}
    |A_{\alpha,\beta}|=\sum \limits_{z\in \mathbb{F}_q^d}\lambda_{r,\theta}(z)^k.
\end{equation}

Plugging $\eqref{eq6.5}$ into inequality $\eqref{eq6.4}$, we obtain 
\begin{equation}
\label{eq6.6}
    |\mathcal{N}_{\theta}(r)|\geq |\Lambda_{\theta}(r)|-\tbinom{k+1}{2}\sum \limits_{z\in \mathbb{F}_q^d}\lambda_{r,\theta}(z)^k.
\end{equation}

Summing $\eqref{eq6.6}$ over all $\theta\in \mathrm{O}_d(\mathbb{F}_q)$, we obtain the following inequality:
\medskip
\begin{equation}
\label{eq6.7}
    \sum \limits_{\theta}|\mathcal{N}_{\theta}(r)|\geq \sum \limits_{\theta}|\Lambda_{\theta}(r)|-\tbinom{k+1}{2}\sum \limits_{\theta,z}\lambda_{r,\theta}(z)^k.
\end{equation}

Combining \eqref{eq6.7} and \eqref{secondmaininequality}, we have
\begin{equation}
\label{finalinequality}
    |\mathcal{S}_A|\geq \frac{1}{|\mathrm{O}_d(\mathbb{F}_q)|}\bigg(\lVert{\lambda_{r,\theta}(z)\rVert}_{k+1}^{k+1}-\tbinom{k+1}{2}\lVert{\lambda_{r,\theta}(z)\rVert}_{k}^{k}\bigg).
\end{equation}

\smallskip

It is worth mentioning that the lower bound for the size of $\mathcal{S}_A$ does not depend on $A$ at all which is the crux of the result being proved.

\smallskip

Now we will show that if $\lvert E\rvert\gg_k q^{d/2}$, then $\lvert \mathcal{S}_{A}\rvert>0$.

\smallskip

Applying Hölder's inequality we obtain the lower bound for the $L^{k+1}$-norm of $\lambda_{r,\theta}(z)$ which is given in the following inequality:
\begin{equation}
\label{holder}
    \sum_{\theta,z}\lambda_{r,\theta}(z)^{k+1}\geq \lvert \mathrm{O}_d(\mathbb{F}_q)\rvert\cdot \frac{|E|^{2k+2}}{q^{dk}}.
\end{equation}

The second summand in \eqref{finalinequality} can be partitioned into two sums where terms are either small or large. More precisely,
\begin{equation}
\label{largesmallparts}
\begin{split}
    \lVert{\lambda_{r,\theta}(z)\rVert}_{k}^{k}&=\sum_{\theta,z}\lambda_{r,\theta}(z)^k \\
    &=\sum_{\lambda\geq c}\lambda_{r,\theta}(z)^k+\sum_{\lambda<c}\lambda_{r,\theta}(z)^k,
\end{split}
\end{equation}
where parameter $c\in \mathbb{R}_{\geq 0}$ will be chosen later. Combining equality \eqref{largesmallparts} with \eqref{finalinequality}, we obtain the following lower bound:
\begin{equation}
\label{lowerboundforSA}
\lvert\mathcal{S}_A\rvert\geq \frac{1}{\lvert \mathrm{O}_d(\mathbb{F}_q)\rvert} \big(\mathcal{S}_1+\mathcal{S}_2\big),    
\end{equation}
where 
\begin{equation*}
    \mathcal{S}_1=\mfrac{1}{2}\lVert \lambda_{r,\theta}(z)\rVert_{k+1}^{k+1}-\tbinom{k+1}{2}\sum_{\lambda\geq c}\lambda_{r,\theta}(z)^k,
\end{equation*}
\smallskip
\begin{equation*}
     \mathcal{S}_2=\mfrac{1}{2}\lVert \lambda_{r,\theta}(z)\rVert_{k+1}^{k+1}-\tbinom{k+1}{2}\sum_{\lambda< c}\lambda_{r,\theta}(z)^k.
\end{equation*}

By definition of the $L^{k+1}$-norm of $\lambda_{r,\theta}(z)$, we have the following lower bound for $\mathcal{S}_1$:
\begin{equation*}
\begin{split}
    \mathcal{S}_1&\geq \mfrac{1}{2}\bigg(\sum_{\lambda\geq c}\lambda_{r,\theta}(z)^{k+1}-k(k+1)\sum_{\lambda\geq c}\lambda_{r,\theta}(z)^{k}\bigg)\\
    &=\mfrac{1}{2}\bigg(\sum_{\lambda \geq c} \lambda_{r,\theta}(z)^k\Big(\lambda_{r,\theta}(z)-k(k+1)\Big)\bigg).
\end{split}
\end{equation*}

Taking parameter $c$ to be equal to $k(k+1)$ it follows that $\mathcal{S}_1\geq 0$.

\smallskip

Now we shall estimate $\mathcal{S}_2$ as follows: for the first term we'll use inequality \eqref{holder} and the second term we'll estimate trivially. Thus,
\begin{equation}
\label{estimateforS2}
\begin{split}
    \mathcal{S}_2&\geq \frac{|E|^{2k+2}}{2q^{dk}}\cdot \lvert \mathrm{O}_d(\mathbb{F}_q)\rvert-\tbinom{k+1}{2}\sum_{\lambda<c}\lambda_{r,\theta}(z)^k \\
    & \geq \frac{|E|^{2k+2}}{2q^{dk}}\cdot \lvert \mathrm{O}_d(\mathbb{F}_q)\rvert - c^k\tbinom{k+1}{2}\sum_{\theta, z} 1\\
    & \geq \lvert \mathrm{O}_d(\mathbb{F}_q)\rvert \Bigg(\frac{|E|^{2k+2}}{2q^{dk}}-\frac{(k^2+k)^{k+1}q^d}{2}\Bigg).
\end{split}
\end{equation}

We notice that $\lvert \mathrm{O}_d(\mathbb{F}_q)\rvert>0$ and hence combining the  inequalities $\mathcal{S}_1\geq 0$ with \eqref{estimateforS2} and \eqref{lowerboundforSA}, we have:
\begin{equation}
    \lvert\mathcal{S}_A\rvert\geq \frac{|E|^{2k+2}}{2q^{dk}}-\frac{(k^2+k)^{k+1}q^d}{2}.
\end{equation}

It is not difficult to verify that if $\lvert E\rvert\geq 2k q^{d/2}$, then $\lvert\mathcal{S}_A\rvert>0$ which completes the proof of Theorem \ref{mainthm}.

\bigskip

\section{The Second Proof of Theorem \ref{mainthm}}
\label{sec:6}

In this section, we'll provide an elementary proof of Theorem \ref{mainthm} which is based on the averaging argument. We'll split the proof into two parts: $r=1$ and $r\in (\mathbb{F}_q)^2\setminus \{0,1\}$.

\smallskip

\underline{\textbf{The case when $r=1$.}} Let $k\geq 1$ and assume that $d\geq 2$ is a positive integer and consider $E\subset \mathbb{F}_q^d$ such that $\lvert E\rvert \geq (k+3) q^{d/2}$. For each $t\in \mathbb{F}_q^d$ we define the translation of set $E$ by an element $t$ as $E+t\coloneqq \big\{x+t: x\in E \big\}$. Assume that $S\subset \mathbb{F}_q^d$ such that $|S|=\lfloor q^{d/2} \rfloor$. 

\smallskip

For each $a\in S$, we have $E+a\subset \mathbb{F}_q^d$. Therefore, we obtain the following set containment: 
\begin{equation}
\label{containemnt}
    \bigcup_{a\in S}(E+a)\subset \mathbb{F}_q^d.
\end{equation}

Comparing cardinality of both sets in \eqref{containemnt} and applying Bonferroni inequality to the RHS, we obtain
\begin{equation*}
    q^d\geq \sum_{a\in S}\lvert E+a\rvert-\sum_{1\leq i<j\leq \lvert S\rvert} \big \lvert(E+a_i)\cap (E+a_j) \big \rvert,
\end{equation*}
where $S=\{a_j\}_{j=1}^{\lvert S\rvert}$.
Since $\lvert E+t\rvert=\lvert E\rvert$, then we have
\begin{equation}
    q^d\geq \lvert S\rvert\lvert E\rvert-\binom{\lvert S\rvert}{2}\max_{1\leq i<j\leq \lvert S\rvert} \big \lvert(E+a_i)\cap (E+a_j) \big \rvert.
\end{equation}

Thus, 
\begin{equation*}
    \max_{1\leq i<j\leq \lvert S\rvert} \big \lvert(E+a_i)\cap (E+a_j) \big \rvert\geq \frac{\lvert S\rvert\lvert E\rvert-q^d}{\tbinom{\lvert S\rvert}{2}}.
\end{equation*}

Using the fact that $\lvert E\rvert\geq (k+3) q^{d/2}$, we have
\begin{equation*}
    \max_{1\leq i<j\leq \lvert S\rvert} \big \lvert(E+a_i)\cap (E+a_j) \big \rvert\geq 2\cdot \frac{(k+3)\lfloor q^{d/2}\rfloor q^{d/2}-q^d}{\lfloor q^{d/2}\rfloor(\lfloor q^{d/2}\rfloor-1)}.
\end{equation*}

Since $x-1<\lfloor x\rfloor\leq x$ and $\lfloor x\rfloor\geq \frac{x}{2}$ for $x\geq 2$, then it follows that 
\begin{equation*}
    \max_{1\leq i<j\leq \lvert S\rvert} \big \lvert(E+a_i)\cap (E+a_j) \big \rvert\geq 2\cdot \frac{(k+3) \cdot \frac{q^{d/2}}{2} \cdot q^{d/2}-q^d}{q^{d/2}(q^{d/2}-1)}
\end{equation*}
\begin{equation*}
    =(k+1)\cdot\frac{q^d}{q^d-q^{d/2}}\geq k+1.
\end{equation*}

\medskip

Therefore, there exists $(i,j)$ such that $1\leq i<j\leq \lvert S\rvert$ such that $ \lvert(E+a_i)\cap (E+a_j)\rvert\geq k+1$. 

\smallskip

We notice that $(E+a_i)\cap (E+a_j)=a_j+E\cap ((a_i-a_j)+E)$ and therefore
\begin{equation*}
    \big\lvert E\cap ((a_i-a_j)+E)\big\rvert\geq k+1.
\end{equation*}

Let $c\coloneqq a_i-a_j\neq 0$ and hence $\lvert E\cap (c+E)\rvert\geq k+1$. It immediately implies that there exists $\{x_i\}_{i=1}^{k+1}$ such that 
\begin{equation*}
    \{x_i\}_{i=1}^{k+1}\subset E\cap (c+E).    
\end{equation*}

Define $y_i\coloneqq x_i-c$ and hence $\{y_i\}_{i=1}^{k+1} \subset E$. Therefore, for each $1\leq i<j \leq k+1$, we have $x_i\neq x_j$, $y_i\neq y_j$ and $\lVert x_i-x_j\rVert=\lVert y_i-y_j\rVert$. Moreover, we notice that $x_i\neq y_i$ for each $i\in [k+1]$ since $c\neq 0$. 

\medskip

\underline{\textbf{The case when $r\in (\mathbb{F}_q)^2\setminus \{0,1\}$.}} So suppose that $r\in (\mathbb{F}_q)^2\setminus \{0,1\}$ and thus $r=t^2$. For this $t$, we define $tE\coloneqq \{tv: v\in E\}$.

Let us consider the following set: 

\begin{equation*}
    H\coloneqq \Big\{(x,a):x\in tE\cap (E+a),\ a\in \mathbb{F}_q^d\Big\}.    
\end{equation*}

We'll utilize the double counting to compute the cardinality of $H$.

\smallskip

Firstly, we see that 
\begin{equation}
\label{sizeofH1}
\begin{split}
    \lvert H\rvert&=\sum_{x\in tE}\sum_{\substack{a\in \mathbb{F}_q^d \\ x\in E+a}}1=\sum_{x\in tE}\sum_{\substack{a\in \mathbb{F}_q^d \\ a\in x-E}}1\\
    &=\sum_{x\in tE}\lvert x-E \rvert=\sum_{x\in tE}\lvert E\rvert\\
    &=\lvert tE\rvert \lvert E\rvert=\lvert E\rvert^2.
\end{split}
\end{equation}

On the other hand, changing the order of variables, we have
\begin{equation}
\label{sizeofH2}
\begin{split}
    \lvert H\rvert&=\sum_{a\in \mathbb{F}_q^d}\sum_{\substack{x\in tE \\ x\in E+a}}1\\
    &=\sum_{a\in \mathbb{F}_q^d}\lvert tE\cap (E+a)\rvert.
\end{split}
\end{equation}

Comparing \eqref{sizeofH1} with \eqref{sizeofH2}, we have
\begin{equation*}
    \sum_{a\in \mathbb{F}_q^d}\lvert tE\cap (E+a)\rvert=\lvert E\rvert^2.
\end{equation*}

Therefore, 
\begin{equation*}
    \max_{a\in \mathbb{F}_q^d}\lvert tE\cap (E+a)\rvert\geq \frac{\lvert E\rvert^2}{q^d}.
\end{equation*}

If $\lvert E\rvert\geq \sqrt{k+2}q^{d/2}$, then 
\begin{equation*}
    \max\limits_{a\in \mathbb{F}_q^d}\lvert tE\cap (E+a)\rvert \geq k+2.    
\end{equation*}

Therefore, there exists $a\in \mathbb{F}_q^d$ such that $\lvert tE\cap (E+a)\rvert \geq k+2$. Thus, there is a sequence $\{x_i\}_{i=1}^{k+2}$ such that $\{x_i\}_{i=1}^{k+2}\subset tE\cap (E+a)$. It implies that $x_i=tz_i$ and $x_i=y_i+a$ for all $i\in [k+2]$, where $y_i,z_i\in E$. In summary, we were able to prove that there are $(y_1,\dots,y_{k+2})\in E^{k+2}$ and $(z_1,\dots,z_{k+2})\in E^{k+2}$ such that $y_i\neq y_j$, $z_i\neq z_j$ if $1\leq i<j\leq k+2$, where $y_i=tz_i-a$ for all $i\in [k+2]$.

\smallskip

Moreover, we notice that there at least $k+1$ $i$'s in $\{1,\dots,k+2\}$ such that $y_i\neq z_i$. Indeed, WLOG assume that $y_1=z_1$. Then for any $j\neq 1$, we have $y_1-y_j=t(z_1-z_j)$ which can be written as $y_1-y_j=t(y_1-z_j)$. It immediately implies that for all $j\neq 1$ we have $y_j\neq z_j$ since by assumption $t\neq 1$.

\bigskip

\section{Corollaries}
\label{sec:4}

In this section, we will describe some interesting corollaries of the Theorem \ref{mainthm}. For brevity, we'll assume that we are in two dimensions and $q$ is a prime number.

\medskip

\textit{\textbf{Paths of length $k$.}} Let $k\geq 2$ and fix any element $r\in (\mathbb{F}_p)^2\setminus\{0\}$. Taking $A=\big\{(1,2),\dots,(k,k+1)\big\}$ we see that if $E\subset \mathbb{F}_p^2$ with $\lvert E\rvert\geq 2kp$, then we can find two $k$-paths in $E$, i.e., $(x_1,\dots,x_{k+1})\in E^{k+1}$ and $(y_1,\dots,y_{k+1})\in E^{k+1}$ such that $x_i\neq x_j,\ y_i\neq y_j$ if $1\leq i<j\leq k+1$ and $\lVert y_i-y_{i+1}\rVert=r\lVert x_i-x_{i+1}\rVert$ if $i\in [k]$.

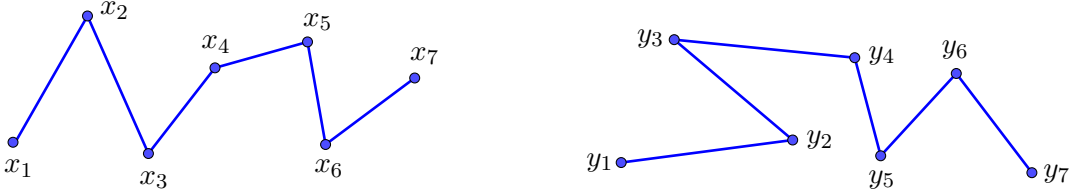
\begin{figure}[htbp]
\centering
\usetikzlibrary{arrows}
\pagestyle{empty}
\definecolor{qqqqff}{rgb}{0.,0.,1.}
\definecolor{ududff}{rgb}{0.30196078431372547,0.30196078431372547,1.}
\begin{tikzpicture}[scale=0.75][line cap=round,line join=round,>=triangle 45,x=1.0cm,y=1.0cm]
\clip(1,-7.5) rectangle (22.5,-4.);
\draw [line width=1.pt,color=qqqqff,fill=qqqqff,fill opacity=0.0] (2.24,-6.64)-- (3.52,-4.4)-- (4.6,-6.84)-- (5.78,-5.32)-- (7.42,-4.86)-- (7.74,-6.68)-- (9.32,-5.5);
\draw [line width=1.pt,color=qqqqff,fill=qqqqff,fill opacity=0.0] (12.98,-7.)-- (16.02,-6.6)-- (13.92,-4.82)-- (17.12,-5.14)-- (17.58,-6.88)-- (18.92,-5.42)-- (20.26,-7.18);
\draw [fill=ududff] (2.2,-6.64) circle (2.5pt);
\draw[color=black] (2.3,-7.1) node {\scalebox{1}{$x_1$}};
\draw [fill=ududff] (3.52,-4.4) circle (2.5pt);
\draw[color=black] (4,-4.3) node {\scalebox{1}{$x_2$}};
\draw [fill=ududff] (4.6,-6.84) circle (2.5pt);
\draw[color=black] (4.7,-7.3) node {\scalebox{1}{$x_3$}};
\draw [fill=ududff] (5.78,-5.32) circle (2.5pt);
\draw[color=black] (5.8,-4.9) node {\scalebox{1}{$x_4$}};
\draw [fill=ududff] (7.42,-4.86) circle (2.5pt);
\draw[color=black] (7.6,-4.5) node {\scalebox{1}{$x_5$}};
\draw [fill=ududff] (7.74,-6.68) circle (2.5pt);
\draw[color=black] (7.8,-7.1) node {\scalebox{1}{$x_6$}};
\draw [fill=ududff] (9.32,-5.5) circle (2.5pt);
\draw[color=black] (9.5,-5.1) node {\scalebox{1}{$x_7$}};
\draw [fill=ududff] (12.98,-7.) circle (2.5pt);
\draw[color=black] (12.6,-7) node {\scalebox{1}{$y_1$}};
\draw [fill=ududff] (16.02,-6.6) circle (2.5pt);
\draw[color=black] (16.5,-6.6) node {\scalebox{1}{$y_2$}};
\draw [fill=ududff] (13.92,-4.82) circle (2.5pt);
\draw[color=black] (13.5,-4.8) node {\scalebox{1}{$y_3$}};
\draw [fill=ududff] (17.12,-5.14) circle (2.5pt);
\draw[color=black] (17.6,-5.1) node {\scalebox{1}{$y_4$}};
\draw [fill=ududff] (17.58,-6.88) circle (2.5pt);
\draw[color=black] (17.6,-7.3) node {\scalebox{1}{$y_5$}};
\draw [fill=ududff] (18.92,-5.42) circle (2.5pt);
\draw[color=black] (18.9,-5) node {\scalebox{1}{$y_6$}};
\draw [fill=ududff] (20.26,-7.18) circle (2.5pt);
\draw[color=black] (20.7,-7.2) node {\scalebox{1}{$y_7$}};
\end{tikzpicture}
\label{6-chains}
\caption{6-paths with dilation ratio $r\in \mathbb{F}_p^{*}$. }
\label{6chainsfigure}
\end{figure}

\medskip

\textit{\textbf{Cycles of length $k$.}} Assume that $k\geq 3$. For $k$-cycles one need to take $A=\big\{(1,2),\dots,(k-1,k), (k,1) \big\}$, then we see that for $E\subset \mathbb{F}_p^2$ with $\lvert E\rvert\geq 2kp$ there are $k$-cycles in $E$, i.e., $(x_1,\dots,x_k)\in E^k$ and $(y_1,\dots,y_k)\in E^k$ such that $x_i\neq x_j,\ y_i\neq y_j$ if $1\leq i<j\leq k$ and $\lVert y_i-y_{i+1}\rVert=r\lVert x_i-x_{i+1}\rVert$ if $i\in \{1,\dots, k-1\}$ and $\lVert y_k-y_{1}\rVert=r\lVert x_k-x_{1}\rVert$.

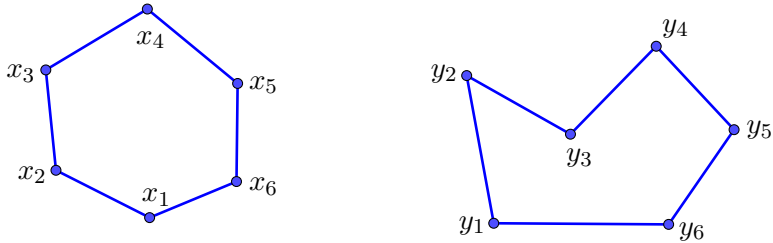
\begin{figure}[htbp]
\centering
\usetikzlibrary{arrows}
\pagestyle{empty}
\definecolor{qqqqff}{rgb}{0.,0.,1.}
\definecolor{ududff}{rgb}{0.30196078431372547,0.30196078431372547,1.}
\begin{tikzpicture}[scale=0.75][line cap=round,line join=round,>=triangle 45,x=1.0cm,y=1.0cm]
\clip(0.,-8.5) rectangle (16.,-4.);
\draw [line width=1.pt,color=qqqqff,fill=qqqqff,fill opacity=0.0] (3.68,-7.96)-- (2.02,-7.12);
\draw [line width=1.pt,color=qqqqff,fill=qqqqff,fill opacity=0.0] (2.02,-7.12)-- (1.84,-5.34);
\draw [line width=1.pt,color=qqqqff,fill=qqqqff,fill opacity=0.0] (1.84,-5.34)-- (3.64,-4.26);
\draw [line width=1.pt,color=qqqqff,fill=qqqqff,fill opacity=0.0] (3.64,-4.26)-- (5.24,-5.58);
\draw [line width=1.pt,color=qqqqff,fill=qqqqff,fill opacity=0.0] (5.24,-5.58)-- (5.22,-7.32);
\draw [line width=1.pt,color=qqqqff,fill=qqqqff,fill opacity=0.0] (5.22,-7.32)-- (3.68,-7.96);
\draw [line width=1.pt,color=qqqqff,fill=qqqqff,fill opacity=0.0] (9.78,-8.06)-- (12.88,-8.08);
\draw [line width=1.pt,color=qqqqff,fill=qqqqff,fill opacity=0.0] (12.88,-8.08)-- (14.04,-6.4);
\draw [line width=1.pt,color=qqqqff,fill=qqqqff,fill opacity=0.0] (14.04,-6.4)-- (12.66,-4.92);
\draw [line width=1.pt,color=qqqqff,fill=qqqqff,fill opacity=0.0] (12.66,-4.92)-- (11.14,-6.48);
\draw [line width=1.pt,color=qqqqff,fill=qqqqff,fill opacity=0.0] (11.14,-6.48)-- (9.3,-5.44);
\draw [line width=1.pt,color=qqqqff,fill=qqqqff,fill opacity=0.0] (9.3,-5.44)-- (9.78,-8.06);
\draw [fill=ududff] (3.68,-7.96) circle (2.5pt);
\draw[color=black] (3.8,-7.6) node {\scalebox{1}{$x_1$}};
\draw [fill=ududff] (2.02,-7.12) circle (2.5pt);
\draw[color=black] (1.6,-7.2) node {\scalebox{1}{$x_2$}};
\draw [fill=ududff] (1.84,-5.34) circle (2.5pt);
\draw[color=black] (1.4,-5.4) node {\scalebox{1}{$x_3$}};
\draw [fill=ududff] (3.64,-4.26) circle (2.5pt);
\draw[color=black] (3.7,-4.8) node {\scalebox{1}{$x_4$}};
\draw [fill=ududff] (5.24,-5.58) circle (2.5pt);
\draw[color=black] (5.7,-5.6) node {\scalebox{1}{$x_5$}};
\draw [fill=ududff] (5.22,-7.32) circle (2.5pt);
\draw[color=black] (5.7,-7.4) node {\scalebox{1}{$x_6$}};
\draw [fill=ududff] (9.78,-8.06) circle (2.5pt);
\draw[color=black] (9.4,-8.1) node {\scalebox{1}{$y_1$}};
\draw [fill=ududff] (12.88,-8.08) circle (2.5pt);
\draw[color=black] (13.3,-8.2) node {\scalebox{1}{$y_6$}};
\draw [fill=ududff] (14.04,-6.4) circle (2.5pt);
\draw[color=black] (14.5,-6.4) node {\scalebox{1}{$y_5$}};
\draw [fill=ududff] (12.66,-4.92) circle (2.5pt);
\draw[color=black] (13,-4.6) node {\scalebox{1}{$y_4$}};
\draw [fill=ududff] (11.14,-6.48) circle (2.5pt);
\draw[color=black] (11.3,-6.9) node {\scalebox{1}{$y_3$}};
\draw [fill=ududff] (9.3,-5.44) circle (2.5pt);
\draw[color=black] (8.9,-5.4) node {\scalebox{1}{$y_2$}};
\end{tikzpicture}
\label{6-cycles}
\caption{6-cycles with dilation ratio $r\in \mathbb{F}_p^{*}$. }
\label{6cyclesfigure}
\end{figure}

\textit{\textbf{Stars with $k$ edges.}} Assume that $k\geq 1$. For $k$-starts one need to take $A=\big\{(1,2),\dots, (1,k+1)\big\}$. Then for $E\subset \mathbb{F}_p^2$ with $|E|\geq 2k p$ there are $k$-stars in $E$, i.e., $(x_1,\dots,x_{k+1})\in E^{k+1}$ and $(y_1,\dots,y_{k+1})\in E^{k+1}$ such that $x_i\neq x_j,\ y_i\neq y_j$ if $1\leq i<j\leq k+1$ and $\lVert y_1-y_i\rVert=r\lVert x_1-x_i\rVert$ if $i\in \{2,\dots,k+1\}$.

\begin{figure}[htbp]
\centering
\usetikzlibrary{arrows}
\pagestyle{empty}
\definecolor{qqqqff}{rgb}{0.,0.,1.}
\definecolor{ududff}{rgb}{0.30196078431372547,0.30196078431372547,1.}
\begin{tikzpicture}[scale=0.75][line cap=round,line join=round,>=triangle 45,x=1.0cm,y=1.0cm]
\clip(-0.5,-6.7) rectangle (16.,-3.);
\draw [line width=1.pt,color=qqqqff,fill=qqqqff,fill opacity=0.0] (3.42,-5.9)-- (0.22,-5.32);
\draw [line width=1.pt,color=qqqqff,fill=qqqqff,fill opacity=0.0] (3.42,-5.9)-- (0.94,-3.94);
\draw [line width=1.pt,color=qqqqff,fill=qqqqff,fill opacity=0.0] (3.42,-5.9)-- (2.98,-3.32);
\draw [line width=1.pt,color=qqqqff,fill=qqqqff,fill opacity=0.0] (3.42,-5.9)-- (4.98,-3.56);
\draw [line width=1.pt,color=qqqqff,fill=qqqqff,fill opacity=0.0] (3.42,-5.9)-- (6.7,-5.54);
\draw [line width=1.pt,color=qqqqff,fill=qqqqff,fill opacity=0.0] (10.942222222222233,-5.857777777777771)-- (14.462222222222236,-5.857777777777771);
\draw [line width=1.pt,color=qqqqff,fill=qqqqff,fill opacity=0.0] (10.942222222222233,-5.857777777777771)-- (14.,-4.968888888888883);
\draw [line width=1.pt,color=qqqqff,fill=qqqqff,fill opacity=0.0] (10.942222222222233,-5.857777777777771)-- (13.53777777777779,-4.231111111111106);
\draw [line width=1.pt,color=qqqqff,fill=qqqqff,fill opacity=0.0] (10.942222222222233,-5.857777777777771)-- (13.05777777777779,-3.8133333333333295);
\draw [line width=1.pt,color=qqqqff,fill=qqqqff,fill opacity=0.0] (10.942222222222233,-5.857777777777771)-- (11.6088888888889,-3.4577777777777743);
\draw [fill=ududff] (3.42,-5.9) circle (2.5pt);
\draw[color=black] (3.5,-6.3) node {\scalebox{1}{$x_1$}};
\draw [fill=ududff] (0.22,-5.32) circle (2.5pt);
\draw[color=black] (0.4,-5) node {\scalebox{1}{$x_2$}};
\draw [fill=ududff] (0.94,-3.94) circle (2.5pt);
\draw[color=black] (1,-3.6) node {\scalebox{1}{$x_3$}};
\draw [fill=ududff] (2.98,-3.32) circle (2.5pt);
\draw[color=black] (3.4,-3.2) node {\scalebox{1}{$x_4$}};
\draw [fill=ududff] (4.98,-3.56) circle (2.5pt);
\draw[color=black] (5.3,-3.3) node {\scalebox{1}{$x_5$}};
\draw [fill=ududff] (6.7,-5.54) circle (2.5pt);
\draw[color=black] (7,-5.3) node {\scalebox{1}{$x_6$}};
\draw [fill=ududff] (10.942222222222233,-5.857777777777771) circle (2.5pt);
\draw[color=black] (10.7,-6.2) node {\scalebox{1}{$y_1$}};
\draw [fill=ududff] (14.462222222222236,-5.857777777777771) circle (2.5pt);
\draw[color=black] (14.9,-5.9) node {\scalebox{1}{$y_6$}};
\draw [fill=ududff] (14.,-4.968888888888883) circle (2.5pt);
\draw[color=black] (14.5,-5) node {\scalebox{1}{$y_5$}};
\draw [fill=ududff] (13.53777777777779,-4.231111111111106) circle (2.5pt);
\draw[color=black] (14,-4.2) node {\scalebox{1}{$y_4$}};
\draw [fill=ududff] (13.05777777777779,-3.8133333333333295) circle (2.5pt);
\draw[color=black] (13.1,-3.4) node {\scalebox{1}{$y_3$}};
\draw [fill=ududff] (11.6088888888889,-3.4577777777777743) circle (2.5pt);
\draw[color=black] (12,-3.4) node {\scalebox{1}{$y_2$}};
\end{tikzpicture}
\label{6-stars}
\caption{5-stars with dilation ratio $r\in \mathbb{F}_p^{*}$. }
\label{5starsfigure}
\end{figure}
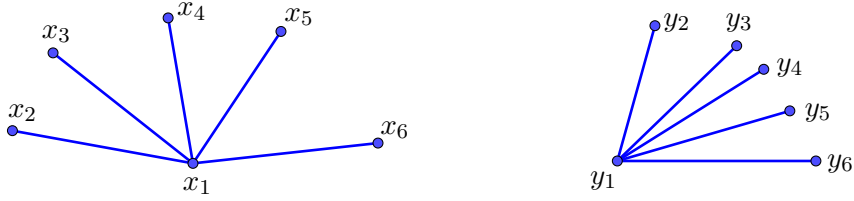

\bigskip

\section{Sharpness example}
\label{sec:5}

In this section we will show that in two dimensions the exponent $d/2$ in Theorem \ref{mainthm} is sharp under some mild conditions.

\smallskip

\textit{\textbf{The case when $\lvert A\rvert\geq 1$.}} We will show that in that setting and if  $q=p^{2\ell}$ with $p\equiv 3 \pmod 4,\ \ell\equiv 1 \pmod 2$ the exponent $d/2$ is sharp in Theorem \ref{mainthm}.

\smallskip

In this setting, our ambient space is $\mathbb{F}_{p^{2\ell}}^{2}$ which is a $2$-dimensional vector space over the finite field $\mathbb{F}_{p^{2\ell}}$ with $p^{2\ell}$ elements. We know that $\mathbb{F}_{p^{2\ell}}$ has the subfield $\mathbb{F}_{p^{\ell}}$ with $p^{\ell}$ elements. Consider the subset $E\coloneqq \mathbb{F}_{p^{\ell}}\times \mathbb{F}_{p^{\ell}}$ of $\mathbb{F}_{p^{2\ell}}^{2}$ and we see that $\lvert E\rvert=q$. 

\smallskip

We know that $\lvert (\mathbb{F}_{p^{2\ell}})^2\setminus \{0\}\rvert=\mfrac{p^{2\ell}-1}{2}>p^{\ell}=\lvert \mathbb{F}_{p^{\ell}}\rvert$ and one can always take $r\in (\mathbb{F}_{p^{2\ell}})^2$ such that $r\notin \mathbb{F}_{p^{\ell}}$. Then for any $(x_1,\dots,x_{k+1})\in E^{k+1},\ (y_1,\dots,y_{k+1})\in E^{k+1}$ and $(\alpha,\beta)\in A$ with $x_i\neq x_j,\ y_i\neq y_j$ for $1\leq i<j\leq k+1$, we have
\begin{equation}
\label{ratio}
    \frac{\lVert y_{\alpha}-y_{\beta}\rVert}{\lVert x_{\alpha}-x_{\beta}\rVert}=r.
\end{equation}

Here we used the fact that $\lVert x\rVert=0$ if and only if $x=(0,0)$ which is true since $p\equiv 3 \pmod 4$ and $\ell\equiv 1 \pmod {2}$. However, in \eqref{ratio} the ratio in the LHS is an element of $\mathbb{F}_{p^{\ell}}$ but $r$ was taken such that $r\notin \mathbb{F}_{p^{\ell}}$. This observation proves that exponent $d/2$ is sharp in this setting. \\

\medskip

\textit{\textbf{The case when $\lvert A\rvert\geq 3$.}} We will show that in that case if $q\equiv 3 \pmod{4}$, then the exponent $d/2$ is sharp in Theorem \ref{mainthm}.

\smallskip

WLOG assume that $(1,2),(1,3),(2,3)\in A$ and assume that $r\in \mathbb{F}_q$ is a nonzero square element such that $r\neq 1$ and $E\coloneqq \big \{(x_1,x_2)\in \mathbb{F}_q^2: x_1^2+x_2^2=1\big\}$ is a sphere of radius $1$ in $\mathbb{F}_q^2$. One can show that $\lvert E\rvert=q+1$. Suppose that there exist $(x_1,x_2,x_3)\in E^3$ and $(y_1,y_2,y_3)\in E^3$ such that $x_i\neq x_j, y_i\neq y_j$ if $1\leq i<j\leq 3$ and $\lVert y_i-y_j\rVert=r\lVert x_i-x_j\rVert$ if $1\leq i<j\leq 3$. For each $i\in [3]$, let's define $z_i\coloneqq \sqrt{r}x_i$ and hence $z_i-z_j=\sqrt{r}(x_i-x_j)$. Therefore, $\lVert y_i-y_j\rVert=\lVert z_i-z_j\rVert$ if $1\leq i<j\leq 3$. It follows that there is $T\in \textup{ISO}(2)$ such that $T(y_i)=z_i$ for $i\in [3]$, where $\textup{ISO}(2)$ is the group of rigid motions.

\smallskip

Since $z_i=\sqrt{r}x_i$, then $\lVert z_i \rVert=r\lVert x_i\rVert=r$. Also we notice that $\lVert z_i-T(\vec{0})\rVert=\lVert y_i\rVert=1$. Therefore, we were able to show that $\{z_j\}_{j=1}^{3}\subset S(T(\vec{0});1)\cap S(\vec{0};r)$, where $S(\vec{a};R)$ denotes the sphere of radius $R$ centered at $\vec{a}$, but we notice that spheres $S(T(\vec{0});1)$ and $S(\vec{0};r)$ are distinct since $r\neq 1$. However, this assumption leads to a contradiction since two distinct spheres in $\mathbb{F}_q^2$ can intersect in at most two points. It shows that the exponent $d/2$ is sharp in this setting.

\bigskip

\bibliographystyle{plain}
\bibliography{refref.bib}

\end{document}